\documentclass[2p]{article}
\usepackage{amssymb,amsmath,amsthm}
\usepackage{indentfirst}
\usepackage{exscale}
\usepackage{relsize}
\usepackage[numbers,sort&compress]{natbib}

\usepackage{geometry}
\usepackage{color}
\newcommand{\R}{\mathbb{R}}
\theoremstyle{definition}

\allowdisplaybreaks
\theoremstyle{remark}

\numberwithin{equation}{section}
\UseRawInputEncoding

\begin{document}
\title{\Large\bf{ Infinitely many solutions for two generalized poly-Laplacian systems on weighted graphs}}
\date{}
\author {Zhangyi Yu$^{1}$, \ Junping Xie$^{2}$\footnote{Corresponding author, E-mail address: hnxiejunping@163.com} , \ Xingyong Zhang$^{1,3}$, \ Wanting Qi$^{1}$\\
{\footnotesize $^1$Faculty of Science, Kunming University of Science and Technology,}\\
 {\footnotesize Kunming, Yunnan, 650500, P.R. China.}\\
  {\footnotesize $^2$Faculty of Transportation Engineering, Kunming University of Science and Technology, Kunming, Yunnan, 650500, P.R. China.}\\
{\footnotesize $^{3}$Research Center for Mathematics and Interdisciplinary Sciences, Kunming University of Science and Technology,}\\
 {\footnotesize Kunming, Yunnan, 650500, P.R. China.}\\}
 \date{}
 \maketitle

 \begin{center}
 \begin{minipage}{15cm}
 \par
 \small  {\bf Abstract:} We investigate the multiplicity of solutions for a generalized poly-Laplacian system on weighted finite graphs and a generalized poly-Laplacian system with Dirichlet boundary value on weighted locally finite graphs, respectively, via the variational methods which are based on mountain pass theorem and topological degree theory. We obtain that these two systems have a sequence of minimax type solutions $\{(u_n,v_n)\}$ satisfying the energy functional $\varphi(u_n,v_n)\to +\infty$ as $n\to +\infty$ and a sequence of local minimum type solutions $\{(u_m^*,v_m^*)\}$ satisfying the energy functional $\varphi(u_m^*,v_m^*)\to -\infty$ as $m\to +\infty$.

 \par
 {\bf Keywords:} generalized poly-Laplacian system, finite graphs, locally finite graphs, mountain pass theorem,  variational methods,  infinitely many solutions.

 \end{minipage}
 \end{center}
  \allowdisplaybreaks
 \vskip2mm
 {\section{Introduction }}
\setcounter{equation}{0}
Assume that $V$ is a vertex set and $E$  is an edge set. Let $G:=(V,E)$ stand for a  graph. If for any
assigned $x \in V$, there are only finite edges $xy \in E$, then $G$ is referred to as a locally finite graph, and if both $V$ and $E$ are finite sets, then $G$ is referred to as a  finite graph. In addition, if any two vertices $x,y\in V$ can be connected by finite edges, then $(V,E)$ is referred to as a connected graph.
\par
In recent years, the problems on the existence of solutions for elliptic equations or systems on finite graph or locally finite graph have attracted some attentions, for example,  \cite{Han2020,Han2021,Yamabe 2016,Han X L2021, Shao M 2023}.  In \cite{Yamabe 2016}, Grigor'Yan-Lin-Yang obtained the existence of a nontrivial solution for the Yamabe type equation on locally finite graph or finite graph by using the mountain pass theorem if the nonlinear term satisfies the super-quadratic condition near $\infty$. Furthermore, they also investigated the $p$-Laplacian equation and poly-Laplacian equation and obtained some similar results. In \cite{Han2021}, Han-Shao studied the $p$-Laplacian equation with Dirichlet boundary condition. They proved that the equation admits a positive solution and a ground state solution by using the mountain pass theorem and the method of Nehari manifold. In \cite{Han2020}, Han-Shao-Zhao studied a nonlinear biharmonic equation with a parameter. When the parameter is small enough, they obtained the existence and convergence of ground state solutions by using the same methods as in \cite{Han2021}. In \cite{Shao M 2023}, Shao investigated the existence of multiple solutions to the nonlinear $p$-Laplacian equation on locally finite graph. Shao proved that the $p$-Laplacian equation admits at least two nontrivial solutions by using the variational methods. In \cite{Han X L2021}, Han-Shao investigated a class of nonlinear $p$-Laplacian equation on a locally finite graph with $p\geq2$. Under appropriate conditions on the nonlinear terms and the potential term, they obtained the equation has a positive solution via the mountain pass theorem.
\par
There have also been some works for the $(p,q)$-Laplacian system and the poly-Laplacian system on graphs, for example, \cite{Shao 2023,Zhang 2022,Yang 2023,Yu 2023,Yang,Pang 2023,Pang 2024}. In \cite{Pang 2023}, by an abstract three critical points theorem without compactness, Pang-Zhang investigated the existence of three distinct solutions for a poly-Laplacian system with a parameter on finite graphs and a $(p,q)$-Laplacian system with a parameter on locally finite graphs. In \cite{Pang 2024}, Pang-Xie-Zhang investigated a generalized poly-Laplacian system with a parameter on weighed finite graphs, a generalized poly-Laplacian system with a parameter and Dirichlet boundary value on weighted locally finite graphs, and a $(p,q)$-Laplacian system with a parameter on weighted locally finite graphs. By utilizing an abstract critical points theorem without compactness condition, they obtained that these systems have infinitely many nontrivial solutions with unbounded norm when the parameters locate some well-determined range. In \cite{Shao 2023}, Shao investigated a class of $p$-Laplacian systems on locally finite graph. By exploiting the method of Nehari manifold and some new analytical techniques, under suitable assumptions on the potentials and nonlinear terms, Shao proved that the  $p$-Laplacian system admits a ground state solution when the parameter is sufficiently large. In \cite{Yang}, Yang-Zhang investigated the existence and multiplicity of nontrivial solutions for a $(p,q)$-Laplacian coupled system with perturbations and two parameters on locally finite graph. By using the Ekeland's variational principle, they obtained that system has at least one nontrivial solution when the nonlinear term satisfies the sub-$(p,q)$ conditions. Moreover, by using the mountain pass theorem and Ekeland's variational principle, they obtained that system has at least one solution of positive energy and one solution of negative energy when the nonlinear term satisfies the super-$(p,q)$ conditions which is weaker than the well-known Ambrosetti-Rabinowitz condition. In \cite{Zhang 2022}, Zhang-Zhang-Xie-Yu investigated the multiplicity of solutions for the following generalized poly-Laplacian system  on a finite graph $G=(V,E)$:
\begin{eqnarray}
\label{eq1}
 \begin{cases}
  \pounds_{m_1,p}u+h_1(x)|u|^{p-2}u=F_u(x,u,v),\;\;\;\;\hfill x\in V,\\
  \pounds_{m_2,q}v+h_2(x)|v|^{q-2}v=F_v(x,u,v),\;\;\;\;\hfill x\in V,
   \end{cases}
\end{eqnarray}
where  $m_i\geq1,\;i=1,2$ are integers,  $1<p,q<\infty$,  $h_i:V\to \R^+,\;i=1,2$, $F:V\times \R \times \R \to \R$, and $\pounds_{m,p}$ is  defined as follows: for any function $\phi:V\to\R$,
\begin{eqnarray}
\label{eq2}
\int_V(\pounds_{m,p}u)\phi d\mu=
 \begin{cases}
  \int_V|\nabla^m u|^{p-2}\Gamma(\Delta^{\frac{m-1}{2}}u,\Delta^{\frac{m-1}{2}}\phi)d\mu,& \text { if $m$ is odd},\\
  \int_V|\nabla^m u|^{p-2}\Delta^{\frac{m}{2}}u\Delta^{\frac{m}{2}}\phi d\mu,&  \text { if $m$ is even}.
    \end{cases}
\end{eqnarray}
When $F$ satisfies the super-$(p,q)$-linear growth condition, by using the mountain pass theorem and the symmetric mountain pass theorem, they obtained the systems have at least one nontrivial solution and $\mbox{dim }W $ nontrivial solutions, where $W$ is the working space of system (\ref{eq1}), respectively. In \cite{Yang 2023},  by using the mountain pass theorem and Ekeland's variational principle, Yang-Zhang investigated the existence of two nontrivial solutions for system (\ref{eq1}) involving concave-convex nonlinearity and parameters with Dirichlet boundary value condition on the locally finite graph. In \cite{Yu 2023}, Yu-Zhang-Xie-Zhang investigated the existence of nontrivial solutions for the following generalized poly-Laplacian system with Dirichlet boundary condition on a locally finite graph $G=(V,E)$:
\begin{eqnarray}
\label{eqA}
 \begin{cases}
   \pounds_{m_1,p}u=F_u(x,u,v),\;\;x\in\Omega,\\
   \pounds_{m_2,q}v=F_v(x,u,v),\;\;x\in\Omega,\\
   |\nabla^j u|=|\nabla^k v|=0,\;\;\;x\in \partial\Omega,\;\;0\le j\le m_1-1,\;\;0\le k\le m_2-1,\\
 \end{cases}
\end{eqnarray}
where $\Omega\cup\partial\Omega\subset V$ is a bounded domain, $m_i ,i=1,2$ are positive integers, and $p=q>1$. By using the mountain pass theorem, when $F$ satisfies the  asymptotically-$p$-linear conditions, they obtained some sufficient conditions about the existence of a nontrivial solution for system (\ref{eqA}).
\par
When $p=2$, $\pounds_{m,p}=(-\Delta)^mu$ is known as the poly-Laplacian operator of $u$, and when $m=1$, $\pounds_{m,p}=-\Delta_pu$ which is known as $p$-Laplacian operator defined as follows:
$$
\Delta_pu(x)=\frac{1}{2\mu(x)}\sum_{y\thicksim x}(|\nabla u|^{p-2}(y)+|\nabla u|^{p-2}(x))\omega_{xy}(u(y)-u(x)).
$$
\par
In this paper, inspired by \cite{Yu 2023,Zhang 2022} but different from them,  we focus on the existence of infinitely many solutions for system (\ref{eq1}) and system (\ref{eqA}), and we allowed $p,q$  to be unequal. The tools such as the symmetric mountain pass theorem, fountain theorem and dual fountain theorem are often used to obtain infinitely many solutions for differential equations in the Euclidean setting, which usually require  the nonlinearity to have symmetry. In this paper, we apply the method in \cite{Habets1995} to the generalized poly-Laplacian systems (\ref{eq1}) and (\ref{eqA}) on graphs, which is essentially based on mountain pass theorem and topological degree theory.  This method has been used to discuss some ordinary differential equations or systems, for example, Che \cite{Che2012}, Habets \cite{Habets1995,Habets1997,Habets1998}, Lu \cite{Haishen 2005}, Jebelean \cite{Jebelean2015}, Li \cite{Li2015} and Liu \cite{Liu2022}.
\par
Finally, we would like to emphasize that it seems that there are no results on infinitely many solutions for the system (\ref{eq1}) on graph excepted for \cite{Pang 2024}. In \cite{Pang 2024}, Pang-Xie-Zhang used the method developed by Bonanno and Bisci in \cite{Bonanno 2010}. They presented some sufficient conditions that (\ref{eq1}) has infinitely many solutions, which is different from our assumptions in Theorem 3.1 and Theorem 4.1 below. The method we used is also different from  that in \cite{Pang 2024}.  Especially,  by this method, we can obtain a sequence of minimax type solutions  and a sequence of local minimum type solutions, which can not be obtained by the methods in  \cite{Pang 2024}.

\vskip2mm
{\section{Preliminaries}}
\setcounter{equation}{0}
\par
In this section, we will review some concepts and properties on poly-Laplacian and Sobolev spaces on graphs. One can see \cite{Yamabe 2016} for more details.
\par
Suppose that $G=(V, E)$ is a graph. The weight on any edge $xy\in E$ is represented by $\omega_{xy}$, which is expected to satisfy $\omega_{xy}=\omega_{yx}$ and $\omega_{xy}>0$. Furthermore, we set $deg(x)=\sum_{y\thicksim x}\omega_{xy}$ for every fixed $x\in V$. We use $y\thicksim x$ to represent $xy\in E$. The distance  of two vertices $x,y$ is denoted by $d(x,y)$, and is defined as the minimal number of edges that connect $x$ with $y$. Let $\Omega\subset V$. If $d(x,y)$ is uniformly bounded for any $x,y\in\Omega$, then $\Omega$ is called  a bounded domain in $V$. The boundary of $\Omega$ is defined by
$$
\partial\Omega=\{y\in V,\;y\notin\Omega|\;\exists\;x\in\Omega\;\text{such\;that}\;xy\in E\},
$$
and the interior of $\Omega$ is defined by $\Omega^{\circ}=\Omega\backslash\partial\Omega.$
\par
Assume that $\mu:V\rightarrow \R^+$ is a finite measure  and there exists a $\mu_0>0$ such that $\mu(x)\ge \mu_0$.
The Laplacian operator $\Delta$ of $\psi$ is defined as
\begin{eqnarray}
\label{eq5}
\Delta \psi(x)=\frac{1}{\mu(x)}\sum\limits_{y\thicksim x}w_{xy}(\psi(y)-\psi(x)).
\end{eqnarray}
Moreover, the corresponding gradient form has the expression given below
\begin{eqnarray}
\label{eq6}
\Gamma(\psi_1,\psi_2)(x)=\frac{1}{2\mu(x)}\sum\limits_{y\thicksim x}w_{xy}(\psi_1(y)-\psi_1(x))(\psi_2(y)-\psi_2(x)).
\end{eqnarray}
Write $\Gamma(\psi)=\Gamma(\psi,\psi)$. The length of the gradient is denoted by
\begin{eqnarray}
\label{eq7}
|\nabla \psi|(x)=\sqrt{\Gamma(\psi)(x)}=\left(\frac{1}{2\mu(x)}\sum\limits_{y\thicksim x}w_{xy}(\psi(y)-\psi(x))^2\right)^{\frac{1}{2}}.
\end{eqnarray}
$|\nabla^m\psi|$ is used to represent the length of the $m$th-order gradient of $\psi$ and $|\nabla^m\psi|$ is defined by
\begin{eqnarray}
\label{eq8}
|\nabla^m\psi|=
 \begin{cases}
  |\nabla\Delta^{\frac{m-1}{2}}\psi|,& \text {if $m$ is odd,}\\
  |\Delta^{\frac{m}{2}}\psi|,&  \text {if $m$ is even,}
   \end{cases}
\end{eqnarray}
here we define $|\nabla\Delta^{\frac{m-1}{2}}\psi|$ by $(\ref{eq7})$ with substituting $\Delta^{\frac{m-1}{2}}\psi$ for $\psi$. For every fixed function $\psi:V\rightarrow\mathbb{R}$, we set
\begin{eqnarray}
\label{eq9}
\int_V \psi(x) d\mu=\sum\limits_{x\in V}\mu(x)\psi(x),
\end{eqnarray}
and $|V|=\sum\limits_{x\in V}\mu(x)$.
\par
For any $p>1$, define
\begin{eqnarray}
\label{eq10}
\Delta_p\psi(x)=\frac{1}{2\mu(x)}\sum\limits_{y\sim x}\left(|\nabla \psi|^{p-2}(y)+|\nabla \psi|^{p-2}(x)\right)\omega_{xy}(\psi(y)-\psi(x)),
\end{eqnarray}
which is called the $p$-Laplacian operator of $\psi$. For any $\phi\in\mathcal{C}_c(V)$, according to the distributional sense, we define $\Delta_p \psi$ as
\begin{eqnarray}
\label{eq11}
\int_V(\Delta_p \psi)\phi d\mu=-\int_V|\nabla \psi|^{p-2}\Gamma(\psi,\phi)d\mu,
\end{eqnarray}
where $\mathcal{C}_c(V)$ is the set of all real functions with compact support. We also observe that $\pounds_{m,p}$ defined by (\ref{eq2}) is the generalization of  $\Delta_p \psi$.

\par
Assume that $G=(V,E)$ is a finite graph. We define
\begin{eqnarray*}
W^{m,l}(V)=\{\psi|\psi:V\to\R\}
\end{eqnarray*}
with the norm
\begin{eqnarray}
\label{eq12}
\|\psi\|_{W^{m,l}(V)}=\left(\int_V(|\nabla^m \psi(x)|^l+h(x)|\psi(x)|^l)d\mu\right)^\frac{1}{l},
\end{eqnarray}
where $m\geq1$, $l>1$, and $h(x)>0$ for all $x\in V$. Then $W^{m,l}(V)$ is Banach space with finite dimension. Let $1<r<+\infty$, set
$$
L^r(V)=\{\psi|\psi:V\to\R\}
$$
equipped by the norm
\begin{eqnarray}
\label{eq13}
\|\psi\|_{r}=\left(\int_V|\psi(x)|^rd\mu\right)^\frac{1}{r}.
\end{eqnarray}

\par
 Each $u\in W^{m,l}(V)$ can be written as $u(x)=\bar u+\tilde u(x)$ with
$$
\bar u:= \frac{1}{|V|} \int_V u(x) d\mu,\;\;\;\int_V \tilde u(x) d\mu=0.
$$
Then $W^{m,l}(V)=\R\oplus\tilde{W}^{m,l}(V)$, where $\tilde{W}^{m,l}(V)=\{\tilde u\in W^{m,l}(V):\int_V \tilde u(x) d\mu=0\}$.

\vskip2mm
\noindent
{\bf Lemma 2.1.} (\cite{Zhang 2022}) {\it Let $l>1$. For all $\psi\in W^{m,l}(V)$, there exists
$$\|\psi\|_{\infty}\leq d_l\|\psi\|_{W^{m,l}(V)},$$
where $d_l=\left(\frac{1}{\mu_{\min}h_{\min}}\right)^\frac{1}{l}$, where $\mu_{\min}:=\min_{x\in V} \mu(x)$ and  $h_{\min}:=\min_{x\in V} h(x)$.}

 \vskip2mm
\noindent
{\bf Lemma 2.2.}  (\cite{Zhang 2022}) {\it Suppose that $G=(V,E)$ is a finite graph. Let $m$ be any positive integer, and let $p>1$. Then $W^{m,p}(V)\hookrightarrow L^q(V)$ for all $1\leq q\leq+\infty$. In particular, if $1< \theta<+\infty$, then for all $\psi\in W^{m,l}(V)$,
\begin{eqnarray}
\label{eq14}
\|\psi\|_{L^\theta(V)}\leq K_l\|\psi\|_{W^{m,l}(V)},
\end{eqnarray}
where
$$
K_l=\frac{\left(\sum_{x\in V}\mu(x)\right)^{\frac{1}{\theta}}}{\mu_{\min}^{\frac{1}{l}}h_{\min}^{\frac{1}{l}}}.
$$
In addition, $W^{m,l}(V)$ is pre-compact, that is, if $\{\psi_k\}$ is bounded in $W^{m,l}(V)$, then up to a subsequence, there exists $\psi\in W^{m,l}(V)$ such that $\psi_k\rightarrow \psi$ in $W^{m,l}(V)$.}

\vskip4mm
\par
Assume that $G(V,E)$ is a locally finite graph and $\Omega\subset G(V,E)$ is a bounded domain. For any given $l>1$ and integer $m\ge1$, we define
$$
W^{m,l}(\Omega)=\{u|u:\Omega\to\R\}
$$
with the norm
\begin{eqnarray}
\label{eqA5}
\|u\|_{W^{m,l}(\Omega)}=\left(\sum_{k=0}^m\int_{\Omega\cup\partial\Omega}|\nabla^ku|^ld\mu\right)^\frac{1}{l}.
\end{eqnarray}
Let
$$
C_c(\Omega)=\{u:V\rightarrow\R|supp \ u\subset\Omega\;\text{and}\;u(x)=0 \mbox{ for } \forall x\in V\backslash\Omega\}.
$$
For any function $\phi\in C_c(\Omega)$, there holds
\begin{eqnarray}
\label{eqA6}
\int_{\Omega}\Delta_lu\phi d\mu=-\int_{\Omega\cup\partial\Omega}|\nabla u|^{l-2}\Gamma(u,\phi) d\mu.
\end{eqnarray}
For any $1\le r<+\infty$, let $L^r(\Omega)$ be the completion of $C_c(\Omega)$ under the norm
$$
\|u\|_{L^r(\Omega)}=\left(\int_{\Omega}|u(x)|^rd\mu\right)^{\frac{1}{r}}.
$$
Let $W_0^{m,l}(\Omega)$ be the completion of $C_c(\Omega)$ under the norm
$$
\|u\|_{W_0^{m,l}(\Omega)}=\left(\int_{\Omega\cup\partial\Omega}|\nabla^mu(x)|^ld\mu\right)^{\frac{1}{s}},
$$
where $m$ is a positive integer and $l>1$. For any $u\in W_0^{m,s}(\Omega)$, we also define the following norm:
$$
\|u\|_{\infty,\Omega}=\max\limits_{x\in\Omega}|u(x)|.
$$
Obviously, $W_0^{m,l}(\Omega)$ is of finite dimension.

 \vskip2mm
\noindent
{\bf Lemma 2.3.}  (\cite{Yamabe 2016}) {\it Suppose that $G=(V,E)$ is a locally finite graph, $\Omega$ is a bounded domain of $V$ such that $\Omega^{\circ}\neq\emptyset$. Let $m$ be any positive integer, and let $l>1$. Then $W_0^{m,l}(\Omega)$ is embedded in $L^{\theta}(\Omega)$ for all $1\leq \theta\leq+\infty$. Particularly, there exists a constant $C_{m,l(\Omega)}>0$ depending only on $m$, $l$ and $\Omega$ such that
\begin{eqnarray}
\label{eqA7}
\left(\int_{\Omega}|u|^\theta d\mu\right)^{\frac{1}{\theta}}\leq C_{m,l(\Omega)}\left(\int_{\Omega\cup\partial\Omega}|\nabla^m u|^ld\mu\right)^{\frac{1}{l}},
\end{eqnarray}
\begin{eqnarray}
\label{eqA8}
\|u\|_{\infty,\Omega}\le H_l\|u\|_{W_0^{m,l}(\Omega)},
\end{eqnarray}
for all $1\le \theta\le +\infty$ and all $u\in W_0^{m,l}(\Omega)$, where $C_{m,l(\Omega)}=\frac{C}{\mu_{\min,\Omega}}(1+|\sum_{x\in\Omega}\mu(x)|)$ with $C$ satisfying $\|u\|_{L^\theta(\Omega)}\le C\|u\|_{W_0^{m,l}(\Omega)}$, $H_l=\frac{C_{m,l}(\Omega)}{\mu_{\min}^\frac{1}{l}}$ and $\mu_{\min,\Omega}=\min_{x\in\Omega}\mu(x)$. Moreover, $W_0^{m,l}(\Omega)$ is pre-compact, that is, if $\{u_n\}$ is bounded in $W_0^{m,l}(\Omega)$, then up to a subsequence, there exists some $u\in W_0^{m,l}(\Omega)$ such that $u_n\rightarrow u$ in $W_0^{m,l}(\Omega)$.}

\vskip2mm
\par
Let $X$ be a Banach space and  $\varphi \in C^{1}(X,\R)$. The functional $\varphi$ satisfies the Palais-Smale (PS) condition if $\{u_n\}$ admits a convergent subsequence in $X$ whenever $\varphi(u_n)$ is bounded and
$\varphi'(u_n)\rightarrow 0$.

 \vskip2mm
\noindent
{\bf Lemma 2.4.} (Mountain pass theorem  \cite{Rabinowitz 1986}) {\it Let $X$ be a real Banach space and $\psi \in C^{1}(X,\R)$, $\psi(0)=0$
and $\psi$ satisfies the (PS)-condition. Suppose that $\psi$ satisfies the following conditions:\\
(i) there exists a constant $ \rho>0$ such that $ \psi|_{\partial B_{\rho}(0)}> 0 $, where $B_\rho=\{a\in X:\|a\|_X<\rho\}$;\\
(ii) there exists $a\in X\backslash \bar B_{\rho} (0)$ such that $\psi(a)\leq 0 $.\\
Then $\psi$ has a critical value $d$ with
 $$
 d:=\inf_{\gamma\in\Gamma}\max_{t\in[0,1]}\psi(\gamma(t)),
 $$
where}
 $$
 \Gamma:=\{\gamma\in C([0,1],X]):\gamma(0)=0,\gamma(1)=a\}.
 $$

\vskip2mm
{\section{The finite graph case}
  \setcounter{equation}{0}
 In this section, we consider the existence of infinitely many solutions for system (\ref{eq1}). First of all, we make the following assumptions:\\
$(A)$ \;  $h_i:V\to\R^+,i=1,2$ and $F(x,s,t)$ is continuously differentiable in $(s,t)\in \R \times \R$ for all $x\in V$;\\
$(H_1)$ \; $\int_V F(x,0,0)d\mu=0$ and there exist functions $f_i:V\to\R^+,i=1,2,3,4$ with $\frac{pq-q}{p}\|f_2\|_{\infty}+\|f_4\|_{\infty}<\frac{1}{K_q^q}$ and $\|f_1\|_{\infty}+\|f_2\|_{\infty}<\frac{1}{K_p^p},$ such that
\begin{eqnarray}
\begin{cases}
|F_s(x,s,t)|\le f_1(x)|s|^{p-1}+f_2(x)|t|^{\frac{pq-q}{p}},\\
|F_t(x,s,t)|\le f_3(x)|s|^{\frac{pq-p}{q}}+f_4(x)|t|^{q-1},\\
\end{cases}
\end{eqnarray}
for all $(x,s,t)\in V\times\R\times\R$, where $K_q=\frac{\left(\sum_{x\in V}\mu(x)\right)^{\frac{1}{q}}}{\mu_{\min}^{\frac{1}{p}}h_{\min}^{\frac{1}{p}}}$ and $K_p=\frac{\left(\sum_{x\in V}\mu(x)\right)^{\frac{1}{p}}}{\mu_{\min}^{\frac{1}{q}}h_{\min}^{\frac{1}{q}}}$};\\
$(H_2)$ \; there hold
\begin{eqnarray}
\label{eq3}
\limsup\limits_{R\rightarrow +\infty}\inf\limits_{a,b\in \R,|(a,b)|=R}\frac{\int_V F(x,a,b)d\mu}{|a|^p+|b|^q} >\max\left\{\frac{1}{p}\int_Vh_1(x)d\mu,\frac{1}{q}\int_Vh_2(x)d\mu\right\}
\end{eqnarray}
and
\begin{eqnarray}
\label{eq4}
& &\liminf\limits_{r\rightarrow +\infty}\sup\limits_{c,d\in \R,|(c,d)|=r}\frac{\int_V F(x,c,d)d\mu}{|c|^p+|d|^q}<\min\Bigg\{\int_V\left(\frac{1}{p}f_1(x)-\frac{1}{p}f_1(x)2^{p-1}-\frac{q-1}{q}f_3(x)-\frac{1}{p}h_1(x)\right)d\mu,\notag\\
& &\int_V\left(\frac{1}{q}f_4(x)-\frac{p-1}{p}f_2(x)2^{q-1}-\frac{1}{q}f_4(x)2^{q-1}-\frac{1}{q}h_2(x)\right)d\mu\Bigg\}.
\end{eqnarray}

\vskip2mm
\par
We shall work in the space $X:=W^{m_1,p}(V)\times W^{m_2,q}(V)$ with the norm $\|(u,v)\|=\|u\|_{W^{m_1,p}(V)}+\|v\|_{W^{m_2,q}(V)}$. Then $X$ is a Banach space with finite dimension. The functional $\varphi_V:X\to\R$ is defined as
\begin{eqnarray}
\label{eq15}
\varphi_V(u,v)=\frac{1}{p}\int_V(|\nabla^{m_1}u|^p+h_1(x)|u|^p)d\mu+\frac{1}{q}\int_V(|\nabla^{m_2}v|^q+h_2(x)|v|^q)d\mu-\int_V F(x,u,v)d\mu.
\end{eqnarray}
Then the assumption $(A)$ implies that $\varphi_V\in C^1(X,\R)$ and
\begin{eqnarray}
\label{eq16}
\langle\varphi_V'(u,v),(\phi_1,\phi_2)\rangle&=&\int_V\left[\pounds_{m_1,p}u\phi_1+h_1(x)|u|^{p-2}u\phi_1-F_u(x,u,v)\phi_1\right]d\mu\nonumber\\
&&+\int_V\left[\pounds_{m_2,q}v\phi_2+h_2(x)|v|^{q-2}v\phi_2-F_v(x,u,v)\phi_2\right]d\mu
\end{eqnarray}
for any $(u,v),(\phi_1,\phi_2)\in X$. The problem about finding the solutions for system (\ref{eq1}) is attributed to finding the critical points of the functional $\varphi_V$ on $X$ (see \cite{Zhang 2022}). The following is our main result in this section.

\vskip2mm
\noindent
{\bf Theorem 3.1.} {\it Suppose that $(A)$, $(H_1)$ and $(H_2)$ hold. Then\\
$(i)$ \; system  (\ref{eq1}) possesses a sequence of solutions  $\{(u_n,v_n)\}$ such that $(u_n,v_n)$ is a mountain pass type critical point of $\varphi_V$ and $\varphi_V(u_n,v_n)=+\infty$  as $n\to \infty$;\\
$(ii)$ \; system  (\ref{eq1}) possesses  a sequence of $\{(u_m^*,v_m^*)\}$  such that $(u_m^*,v_m^*)$ is a local minimum point of $\varphi_V$ and $\varphi_V(u_m^*,m_m^*)=-\infty$ as $m\to \infty$.}

\vskip2mm
  \par
Theorem 3.1 will be proved from the following lemmas. Consider the direct sum decomposition $X=\R^2\oplus Q_V$ with $Q_V=\{(u,v)\in X:\int_V u(x)d\mu=0,\int_Vv(x) d\mu=0\}$. A simple calculation implies that for all $(u,v)\in X$ we can write $(u,v)=(\bar{u},\bar{v})+(\tilde{u},\tilde{v})$ with $(\bar{u},\bar{v})\in \R^2$ and $(\tilde{u},\tilde{v})\in Q_V$. Moreover, we denote $B_{R}=\{(u,v)\in \R^2: |(u,v)|\le R\}$ and $\partial B_{R}$ denotes the boundary of $B_{R}$.

 \vskip2mm
 \noindent
{\bf Lemma  3.1.}  {\it $\varphi_V$ is coercive on $X$.}
 \vskip0mm
 \noindent
{\bf Proof.}\;Note that $\int_V F(x,0,0)d\mu=0$. Then by $(H_1)$, we have
\begin{eqnarray*}
\int_VF(x,s,t)d\mu&=&\int_VF(x,s,t)d\mu-\int_VF(x,0,0)d\mu\\
&=&\int_V(F(x,s,t)-F(x,0,t))d\mu+\int_V(F(x,0,t)-F(x,0,0))d\mu\\
&\le&\int_V|F(x,s,t)-F(x,0,t)|d\mu+\int_V|F(x,0,t)-F(x,0,0)|d\mu\\
&\le&\int_V\int_0^{|s|}|F_s(x,s,t)|dsd\mu+\int_V\int_0^{|t|}|F_t(x,0,t)|dtd\mu\\
&\le&\int_V\int_0^{|s|}(f_1(x)|s|^{p-1}+f_2(x)|t|^{\frac{pq-q}{p}})dsd\mu+\int_V\int_0^{|t|}(f_4(x)|t|^{q-1})dtd\mu\\
&=&\int_V\left(\frac{1}{p}f_1(x)|s|^p+f_2(x)|t|^{\frac{pq-q}{p}}|s|+\frac{1}{q}f_4(x)|t|^q\right)d\mu\\
&\le&\frac{1}{p}\|f_1\|_{\infty}\|s\|_p^p+\|f_2\|_{\infty}\int_V|t|^{\frac{pq-q}{p}}|s|d\mu+\frac{1}{q}\|f_4\|_{\infty}\|t\|_q^q.\\
\end{eqnarray*}
Thus
\begin{eqnarray*}
\varphi_V(u,v)&=&\frac{1}{p}\|u\|_{W^{m_1,p}(V)}^p+\frac{1}{q}\|v\|_{W^{m_2,q}(V)}^q-\int_V F(x,u,v)d\mu\\
&\ge&\frac{1}{p}\|u\|_{W^{m_1,p}(V)}^p+\frac{1}{q}\|v\|_{W^{m_2,q}(V)}^q-\frac{1}{p}\|f_1\|_{\infty}\|u\|_p^p-\|f_2\|_{\infty}\int_V|v|^{\frac{pq-q}{p}}|u|d\mu-\frac{1}{q}\|f_4\|_{\infty}\|v\|_q^q\\
&\ge&\frac{1}{p}\|u\|_{W^{m_1,p}(V)}^p+\frac{1}{q}\|v\|_{W^{m_2,q}(V)}^q-\frac{1}{p}\|f_1\|_{\infty}\|u\|_p^p-\|f_2\|_{\infty}\int_V\left(\frac{|v|^q}{\frac{p}{p-1}}+\frac{|u|^p}{p}\right)d\mu
-\frac{1}{q}\|f_4\|_{\infty}\|v\|_q^q\\
&=&\frac{1}{p}\|u\|_{W^{m_1,p}(V)}^p+\frac{1}{q}\|v\|_{W^{m_2,q}(V)}^q-\frac{1}{p}\|f_1\|_{\infty}\|u\|_p^p-\frac{p-1}{p}\|f_2\|_{\infty}\|v\|_q^q-\frac{1}{p}\|f_2\|_{\infty}\|u\|_p^p-\frac{1}{q}\|f_4\|_{\infty}\|v\|_q^q\\
&\ge&\frac{1}{p}\|u\|_{W^{m_1,p}(V)}^p+\frac{1}{q}\|v\|_{W^{m_2,q}(V)}^q-\frac{1}{p}\|f_1\|_{\infty}K_p^p\|u\|_{W^{m_1,p}(V)}^p-\frac{p-1}{p}\|f_2\|_{\infty}K_q^q\|v\|_{W^{m_2,q}(V)}^q\\
& &-\frac{1}{p}\|f_2\|_{\infty}K_p^p\|u\|_{W^{m_1,p}(V)}^p-\frac{1}{q}\|f_4\|_{\infty}K_q^q\|v\|_{W^{m_2,q}(V)}^q.
\end{eqnarray*}
Since $\frac{pq-q}{p}\|f_2\|_{\infty}+\|f_4\|_{\infty}<\frac{1}{K_q^q}$ and $\|f_1\|_{\infty}+\|f_2\|_{\infty}<\frac{1}{K_p^p}$, then $\varphi_V(u,v)$ is coercive on $X$.
\qed
 \vskip2mm
 \noindent
{\bf Lemma 3.2.}  {\it There is a sequence $\{R_n\}$ such that
$$
\lim\limits_{n\to\infty}R_n=+\infty,\;\lim\limits_{n\to\infty}\left[\sup\limits_{(a,b)\in\R^2,|(a,b)|=R_n}\varphi_V(a,b)\right]=-\infty.
$$}
 \vskip0mm
 \noindent
{\bf Proof.}\; Note that $|\nabla^{m_1} a|=|\nabla^{m_2} b|=0$ for all $a,b\in\R$. Thus
\begin{eqnarray*}
        \varphi_V(a,b)
&  =  & \frac{|a|^p}{p}\int_Vh_1(x)d\mu+\frac{|b|^q}{q}\int_Vh_2(x)d\mu-\int_V F(x,a,b) d\mu\\
& \le &  \max\left\{\frac{1}{p}\int_Vh_1(x)d\mu,\frac{1}{q}\int_Vh_2(x)d\mu\right\}(|a|^p+|b|^q)-\int_V F(x,a,b) d\mu\\
&  =  & (|a|^p+|b|^q)\left(\max\left\{\frac{1}{p}\int_Vh_1(x)d\mu,\frac{1}{q}\int_Vh_2(x)d\mu\right\}-\frac{\int_V F(x,a,b)d\mu}{|a|^p+|b|^q}\right).
\end{eqnarray*}
Then  the proof is completed by  (\ref{eq3}).\qed

\vskip2mm
 \noindent
{\bf Proposition 3.1.}  {\it For any $c\in\R$, there holds $|\nabla^m(c+\tilde{u})|=|\nabla^m\tilde{u}|$.}
\vskip0mm
 \noindent
{\bf Proof.}\; Firstly, we prove that there holds $\Delta(\Delta^nu)=\Delta^n(\Delta u)$ for any $n\in\mathbb{N}$. Obviously, it holds when $n=1$,
\begin{eqnarray*}
\Delta(\Delta u)=\Delta(\Delta u).
\end{eqnarray*}
Suppose that there holds $\Delta(\Delta^{k-1}u)=\Delta^{k-1}(\Delta u)$. Then when $n=k$,
\begin{eqnarray*}
\Delta(\Delta^ku)=\Delta(\Delta(\Delta^{k-1} u))=\Delta(\Delta^{k-1}(\Delta u))=\Delta^k(\Delta u).
\end{eqnarray*}
Thus $\Delta(\Delta^nu)=\Delta^n(\Delta u)$ holds by mathematical induction.
\par
Secondly, for any $c\in\R$, there hold
\begin{eqnarray*}
\Delta(c+\tilde{u})(y)=\frac{1}{\mu(y)}\sum\limits_{z\thicksim y}w_{yz}(\tilde{u}(z)-\tilde{u}(y))=\Delta\tilde{u}(y),\\
\Delta(c+\tilde{u})(x)=\frac{1}{\mu(x)}\sum\limits_{y\thicksim x}w_{xy}(\tilde{u}(y)-\tilde{u}(x))=\Delta\tilde{u}(x).
\end{eqnarray*}
When $m=1$, it is obvious that
\begin{eqnarray*}
|\nabla(c+\tilde{u})|=|\nabla\tilde{u}|.
\end{eqnarray*}
When $m=2$,
\begin{eqnarray*}
|\nabla^2\tilde{u}|=|\Delta\tilde{u}|=\left|\frac{1}{\mu(x)}\sum\limits_{y\thicksim x}w_{xy}(\tilde{u}(y)-\tilde{u}(x))\right|=|\Delta(c+\tilde{u})|=|\nabla^2(c+\tilde{u})|.
\end{eqnarray*}
When $m=3$,
\begin{eqnarray*}
& &|\nabla^3(c+\tilde{u})|=|\nabla\Delta(c+\tilde{u})|=\left(\frac{1}{2\mu(x)}\sum\limits_{y\thicksim x}w_{xy}(\Delta(c+\tilde{u})(y)-\Delta(c+\tilde{u})(x))^2\right)^{\frac{1}{2}}\\
&=&\left(\frac{1}{2\mu(x)}\sum\limits_{y\thicksim x}w_{xy}(\Delta\tilde{u}(y)-\Delta\tilde{u}(x))^2\right)^{\frac{1}{2}}=|\nabla\Delta\tilde{u}|=|\nabla^3\tilde{u}|.
\end{eqnarray*}
When $m=2k+2\;(k\in\mathbb{N})$,
\begin{eqnarray*}
|\nabla^{2k+2}(c+\tilde{u})|=|\Delta^{k+1}(c+\tilde{u})|=|\Delta(\Delta^k(c+\tilde{u}))|=|\Delta^k(\Delta(c+\tilde{u}))|=|\Delta^k(\Delta\tilde{u})|=|\Delta^{k+1}\tilde{u}|=|\nabla^{2k+2}\tilde{u}|.
\end{eqnarray*}
When $m=2k+3\;(k\in\mathbb{N})$,
\begin{eqnarray*}
|\nabla^{2k+3}(c+\tilde{u})|=|\nabla\Delta^{k+1}(c+\tilde{u})|=|\nabla\Delta^k(\Delta(c+\tilde{u}))|=|\nabla\Delta^k(\Delta\tilde{u})|=|\nabla\Delta^{k+1}\tilde{u}|=|\nabla^{2k+3}\tilde{u}|.
\end{eqnarray*}
Hence, for any $c\in\R$ and $m\in\mathbb{N}$, there holds $|\nabla^m(c+\tilde{u})|=|\nabla^m\tilde{u}|$.
\qed

\vskip2mm
 \noindent
{\bf Lemma 3.3.}  {\it There is a sequence $\{r_m\}$ such that
$$
\lim\limits_{m\to\infty}r_m=+\infty,\;\lim\limits_{m\to\infty}\left[\inf\limits_{(c,d)\in\R^2,|(c,d)|=r_m,(\tilde{u},\tilde{v})\in Q_V}\varphi_V(c+\tilde{u},d+\tilde{v})\right]=+\infty.
$$}
 \vskip0mm
 \noindent
{\bf Proof.}\; For any given $c,d\in\R$, $|(c,d)|=r_m$, and $(\tilde{u},\tilde{v})\in Q_V$, set $(u,v)=(c+\tilde{u},d+\tilde{v})$. We have
\begin{eqnarray*}
& &\int_V F(x,c+\tilde{u},d+\tilde{v})d\mu\\
&=&\int_VF(x,c,d)d\mu+\int_V[F(x,c+\tilde{u},d+\tilde{v})-F(x,c,d)]d\mu\\
&=&\int_VF(x,c,d)d\mu+\int_V\int_0^1F_s(x,c+s\tilde{u},d+\tilde{v})\tilde{u}dsd\mu+\int_V\int_0^1F_t(x,c,d+t\tilde{v})\tilde{v}dtd\mu\\
&\ge&\int_VF(x,c,d)d\mu+\int_V\int_0^1\left(f_1(x)|c+s\tilde{u}|^{p-1}+f_2(x)|d+\tilde{v}|^{\frac{pq-q}{p}}\right)\tilde{u}dsd\mu\\
& &+\int_V\int_0^1\left(f_3(x)|c|^{\frac{pq-p}{q}}+f_4(x)|d+t\tilde{v}|^{q-1}\right)\tilde{v}dtd\mu\\
&=&\int_VF(x,c,d)d\mu+\int_V\left(\frac{1}{p}f_1(x)|c+\tilde{u}|^p-\frac{1}{p}f_1(x)|c|^p+f_2(x)|d+\tilde{v}|^{\frac{pq-q}{p}}\tilde{u}\right)d\mu\\
& &+\int_V\left(\frac{1}{q}f_4(x)|d+\tilde{v}|^q-\frac{1}{q}f_4(x)|d|^q+f_3(x)|c|^{\frac{pq-p}{q}}\tilde{v}\right)d\mu\\
&\ge&\int_VF(x,c,d)d\mu+\int_Vf_2(x)|d+\tilde{v}|^{\frac{pq-q}{p}}|\tilde{u}|d\mu+\int_Vf_3(x)|c|^{\frac{pq-p}{q}}|\tilde{v}|d\mu\\
& &+\int_V\left[\frac{1}{p}f_1(x)2^{p-1}(|c|^p+|\tilde{u}|^p)-\frac{1}{p}f_1(x)|c|^p\right]d\mu+\int_V\left[\frac{1}{q}f_4(x)2^{q-1}(|d|^q+|\tilde{v}|^q)-\frac{1}{q}f_4(x)|d|^q\right]d\mu\\
&\ge&\int_VF(x,c,d)d\mu+\int_Vf_2(x)\left(\frac{|d+\tilde{v}|^q}{\frac{p}{p-1}}+\frac{|\tilde{u}|^p}{p}\right)d\mu+\int_Vf_3(x)\left(\frac{|c|^p}{\frac{q}{q-1}}+\frac{|\tilde{v}|^q}{q}\right)d\mu\\
& &+\int_V\left(\frac{1}{p}f_1(x)2^{p-1}|\tilde{u}|^p+\frac{1}{q}f_4(x)2^{q-1}|\tilde{v}|^q\right)d\mu\\
& &+|c|^p\int_V\left(\frac{1}{p}f_1(x)2^{p-1}-\frac{1}{p}f_1(x)\right)d\mu+|d|^q\int_V\left(\frac{1}{q}f_4(x)2^{q-1}-\frac{1}{q}f_4(x)\right)d\mu\\
&\ge&\int_VF(x,c,d)d\mu+\int_V\bigg[\frac{1}{p}(f_1(x)2^{p-1}+f_2(x))|\tilde{u}|^p+\frac{1}{q}\left(f_3(x)+2^{q-1}f_4(x)+2^{q-1}\frac{pq-q}{p}f_2(x)\right)|\tilde{v}|^q\bigg]d\mu\\
& &+|c|^p\int_V\left(\frac{1}{p}f_1(x)2^{p-1}-\frac{1}{p}f_1(x)+\frac{q-1}{q}f_3(x)\right)d\mu+|d|^q\int_V\left(\frac{p-1}{p}f_2(x)2^{q-1}+\frac{1}{q}f_4(x)2^{q-1}-\frac{1}{q}f_4(x)\right)d\mu\\
&\ge&\int_VF(x,c,d)d\mu+\frac{1}{p}(\|f_1\|_{\infty}2^{p-1}+\|f_2\|_{\infty})K_p^p\|\tilde{u}\|_{W^{m_1,p}(V)}^p\\
& &+\frac{1}{q}\left(
\|f_3\|_{\infty}+2^{q-1}\|f_4\|_{\infty}+2^{q-1}\frac{pq-q}{p}\|f_2\|_{\infty}\right)K_q^q\|\tilde{v}\|_{W^{m_2,q}(V)}^q\\
& &+|c|^p\int_V\left(\frac{1}{p}f_1(x)2^{p-1}-\frac{1}{p}f_1(x)+\frac{q-1}{q}f_3(x)\right)d\mu+|d|^q\int_V\left(\frac{p-1}{p}f_2(x)2^{q-1}+\frac{1}{q}f_4(x)2^{q-1}-\frac{1}{q}f_4(x)\right)d\mu.
\end{eqnarray*}
Thus
\begin{eqnarray*}
& &\varphi_V(u,v)\\
&=&\varphi_V(c+\tilde{u},d+\tilde{v})\\
&=&\frac{1}{p}\int_V(|\nabla^{m_1}(c+\tilde{u})|^p+h_1(x)|c+\tilde{u}|^p)d\mu+\frac{1}{q}\int_V(|\nabla^{m_2}(d+\tilde{v})|^q+h_2(x)|d+\tilde{v}|^q)d\mu-\int_V F(x,c+\tilde{u},d+\tilde{v})d\mu\\
&=&\frac{1}{p}\int_V(|\nabla^{m_1}\tilde{u}|^p+h_1(x)|c+\tilde{u}|^p)d\mu+\frac{1}{q}\int_V(|\nabla^{m_2}\tilde{v}|^q+h_2(x)|d+\tilde{v}|^q)d\mu-\int_V F(x,c+\tilde{u},d+\tilde{v})d\mu\\
&\ge&\frac{1}{p}\int_V(|\nabla^{m_1}\tilde{u}|^p+h_1(x)|c+\tilde{u}|^p)d\mu-\frac{1}{p}(\|f_1\|_{\infty}2^{p-1}+\|f_2\|_{\infty})K_p^p\int_V(|\nabla^{m_1}\tilde{u}|^p+h_1(x)|\tilde{u}|^p)d\mu\\
& &+\frac{1}{q}\int_V(|\nabla^{m_2}\tilde{v}|^q+h_2(x)|d+\tilde{v}|^q)d\mu\\
& &-\frac{1}{q}\left(
\|f_3\|_{\infty}+2^{q-1}\|f_4\|_{\infty}+2^{q-1}\frac{pq-q}{p}\|f_2\|_{\infty}\right)K_q^q\int_V(|\nabla^{m_2}\tilde{v}|^q+h_2(x)|\tilde{v}|^q)d\mu\\
& &-|c|^p\int_V\left(\frac{1}{p}f_1(x)2^{p-1}-\frac{1}{p}f_1(x)+\frac{q-1}{q}f_3(x)\right)d\mu\\
& &-|d|^q\int_V\left(\frac{p-1}{p}f_2(x)2^{q-1}+\frac{1}{q}f_4(x)2^{q-1}-\frac{1}{q}f_4(x)\right)d\mu-\int_V F(x,c,d)d\mu\\
&\ge&\frac{1}{p}\int_V\left(|\nabla^{m_1}\tilde{u}|^p+h_1(x)(2^{1-p}|\tilde{u}|^p-|c|^p)\right)d\mu-\frac{1}{p}(\|f_1\|_{\infty}2^{p-1}+\|f_2\|_{\infty})K_p^p\int_V(|\nabla^{m_1}\tilde{u}|^p+h_1(x)|\tilde{u}|^p)d\mu\\
& &+\frac{1}{q}\int_V\left(|\nabla^{m_2}\tilde{v}|^q+h_2(x)(2^{1-q}|\tilde{v}|^q-|d|^q)\right)d\mu\\
& &-\frac{1}{q}\left(
\|f_3\|_{\infty}+2^{q-1}\|f_4\|_{\infty}+2^{q-1}\frac{pq-q}{p}\|f_2\|_{\infty}\right)K_q^q\int_V(|\nabla^{m_2}\tilde{v}|^q+h_2(x)|\tilde{v}|^q)d\mu\\
& &-|c|^p\int_V\left(\frac{1}{p}f_1(x)2^{p-1}-\frac{1}{p}f_1(x)+\frac{q-1}{q}f_3(x)\right)d\mu\\
& &-|d|^q\int_V\left(\frac{p-1}{p}f_2(x)2^{q-1}+\frac{1}{q}f_4(x)2^{q-1}-\frac{1}{q}f_4(x)\right)d\mu-\int_V F(x,c,d)d\mu\\
&=&\frac{1}{p}[1-(\|f_1\|_{\infty}2^{p-1}+\|f_2\|_{\infty})K_p^p]\int_V|\nabla^{m_1}\tilde{u}|^p d\mu+\frac{1}{p}[2^{1-p}-(\|f_1\|_{\infty}2^{p-1}+\|f_2\|_{\infty})K_p^p]\int_V h_1(x)|\tilde{u}|^p d\mu\\
& &+\frac{1}{q}\bigg[1-\left(\|f_3\|_{\infty}+2^{q-1}\|f_4\|_{\infty}+2^{q-1}\frac{pq-q}{p}\|f_2\|_{\infty}\right)K_q^q\bigg]\int_V|\nabla^{m_2}\tilde{v}|^q d\mu\\
& &+\frac{1}{q}\bigg[2^{1-q}-\left(\|f_3\|_{\infty}+2^{q-1}\|f_4\|_{\infty}+2^{q-1}\frac{pq-q}{p}\|f_2\|_{\infty}\right)K_q^q\bigg]\int_V h_2(x)|\tilde{v}|^q d\mu\\
& &-|c|^p\int_V\left(\frac{1}{p}h_1(x)+\frac{1}{p}f_1(x)2^{p-1}-\frac{1}{p}f_1(x)+\frac{q-1}{q}f_3(x)\right)d\mu\\
& &-|d|^q\int_V\left(\frac{1}{q}h_2(x)+\frac{p-1}{p}f_2(x)2^{q-1}+\frac{1}{q}f_4(x)2^{q-1}-\frac{1}{q}f_4(x)\right)d\mu-\int_V F(x,c,d)d\mu\\
&\ge&\frac{1}{p}[1-(\|f_1\|_{\infty}2^{p-1}+\|f_2\|_{\infty})K_p^p]\int_V|\nabla^{m_1}\tilde{u}|^p d\mu+\frac{1}{p}[2^{1-p}-(\|f_1\|_{\infty}2^{p-1}+\|f_2\|_{\infty})K_p^p]\int_V h_1(x)|\tilde{u}|^p d\mu\\
& &+\frac{1}{q}\bigg[1-\left(\|f_3\|_{\infty}+2^{q-1}\|f_4\|_{\infty}+2^{q-1}\frac{pq-q}{p}\|f_2\|_{\infty}\right)K_q^q\bigg]\int_V|\nabla^{m_2}\tilde{v}|^q d\mu\\
& &+\frac{1}{q}\bigg[2^{1-q}-\left(\|f_3\|_{\infty}+2^{q-1}\|f_4\|_{\infty}+2^{q-1}\frac{pq-q}{p}\|f_2\|_{\infty}\right)K_q^q\bigg]\int_V h_2(x)|\tilde{v}|^q d\mu\\
& &+(|c|^p+|d|^q)\Bigg(\min\bigg\{\int_V\left(\frac{1}{p}f_1(x)-\frac{1}{p}f_1(x)2^{p-1}-\frac{q-1}{q}f_3(x)-\frac{1}{p}h_1(x)\right)d\mu,\\
& &\int_V\left(\frac{1}{q}f_4(x)-\frac{p-1}{p}f_2(x)2^{q-1}-\frac{1}{q}f_4(x)2^{q-1}-\frac{1}{q}h_2(x)\right)d\mu\bigg\}-\frac{\int_VF(x,c,d)d\mu}{|c|^p+|d|^q}\Bigg).
\end{eqnarray*}
Therefore, the proof is completed by  (\ref{eq4}).\qed
\par
 \vskip2mm
Consider the set $S_n=\{\gamma\in C(B_{R_n},X),\gamma|_{\partial B_{R_n}}=id|_{\partial B_{R_n}}\}$, where $R_n$ is defined in Lemma 3.2. Define
\begin{eqnarray}
\label{eq20}
c_n=\inf_{\gamma\in\ S_n}\left[\max_{(x,y)\in B_{R_n}}\varphi_V(\gamma(x,y))\right].
\end{eqnarray}
We claim that each $\gamma$ intersects $Q_V$. Suppose that $\pi_V:X\to \R^2$ is the (continuous) projection of $X$ onto $\R^2$, which is defined by
$$
\pi_V(u,v)=\left(\frac{1}{|V|}\int_V u d\mu,\frac{1}{|V|}\int_V v d\mu\right)\;\;\;\mbox{for}\;(u,v)\in X.
$$
Assume that $\gamma$ is any continuous map such that $\gamma|_{\partial B_{R_n}}=id|_{\partial B_{R_n}}$. We need to prove that $(0,0)\in\pi(\gamma(B_{R_n}))$.
\par
Define
$$
\gamma_t(u,v)=t\pi_V(\gamma(u,v))+(1-t)(u,v),
$$
where $t\in[0,1]$ and $(u,v)\in \R^2$. Note that $\gamma_t \in C^0(\R^2;\R^2)$ defines a homotopy of $\gamma_0=id$ with $\gamma_1=\pi_V\circ\gamma$. Furthermore by  $\gamma|_{\partial B_{R_n}}=id|_{\partial B_{R_n}}$ and $(u,v)\in \R^2$, we can obtain that $\gamma_t|_{\partial B_{R_n}}=id$ for all $t$. Thus, we conclude that
$$
deg(\pi_V\circ\gamma,B_{R_n},(0,0))=deg(id,B_{R_n},(0,0))=1,
$$
by homotopy invariance and normalization of the degree (see for instance Deimling \cite{K Deimling 1985}). Then there exists $p^*\in B_{R_n}$ such that $\pi\circ\gamma(p^*)=(0,0)$. So we have $\gamma(p^*)\in Q_V$ (Otherwise if $\gamma(p^*)\in \R^2$, by $\pi\circ\gamma(p^*)=(0,0)$ and the definitions of $\pi$ and $\R^2$ we also have $\gamma(p^*)=(0,0)\in Q_V$). Consequently, each $\gamma$ intersects $Q_V$, that is $\gamma(B_{R_n})\cap Q_V\not=\emptyset$.
\par
By Lemma 3.1, $\varphi_V$ is coercive on $X$ and of course, also on $Q_V$. Then by the continuity of $\varphi$, there exists a constant $M$ such that
$$
\max_{(x,y)\in B_{R_n}}\varphi_V(\gamma(x,y))\geq\max_{(u,v)\in\gamma(B_{R_n})\cap Q_V}\varphi_V(u,v)\geq\inf_{(\tilde{u},\tilde{v})\in Q_V}\varphi_V(\tilde{u},\tilde{v})\geq M.
$$
Therefore $c_n\geq M$, which together with Lemma 3.2 implies that
\begin{eqnarray}
\label{eq20B}
c_n >\sup_{(a,b)\in\R^2,|(a,b)|=R_n}\varphi_V(a,b),
\end{eqnarray}
for all large values of $n$. The following result can be obtained by utilizing $(\ref{eq20B})$ and Theorem 4.3 in \cite{Mawhin1989}.

\vskip2mm
 \noindent
{\bf Lemma 3.4.}  {\it For any given $n$ large enough, assume that there are sequences $(\gamma_k)$ in $S_n$ and $\{(u_n,v_n)\}$ in $X$ such that
$$
\max_{(x,y)\in B_{R_n}}\varphi_V(\gamma_k(x,y))\to c_n\;\;\; \mbox{as }k\to \infty.
$$
Then there is a sequence $\{(u_k,v_k)\}$ in $X$ such that
$$
\varphi_V(u_k,v_k)\to c_n,\;\;\;dist((u_k,v_k),\gamma_k(B_{R_n}))\to 0,\;\;\;|\varphi_V'(u_k,v_k)|\to 0
$$
when $k\to\infty.$
}

\vskip2mm
 \noindent
{\bf Lemma 3.5.}  {\it The sequence $\{(u_k,v_k)\}$ is bounded in $X$.}
 \vskip0mm
 \noindent
{\bf Proof.}\; By Lemma 3.4, when $k$ is large enough, we have
$$
c_n\le\max_{(x,y)\in B_{R_n}}\varphi_V(\gamma_k(x,y))\le c_n+1
$$
and we can take $(w_k,z_k)\in\gamma_k(B_{R_n})$ such that
\begin{eqnarray}
\label{eq21}
\|(u_k,v_k)-(w_k,z_k)\|_X\le1.
\end{eqnarray}
Since
$$
\varphi_V(w_k,z_k)\le c_n+1
$$
and $\varphi_V$ is coercive on $X$, the sequence $\{(w_k,z_k)\}$ is bounded for all $k\in\mathbb{N}$ in $X$, and then the lemma follows from (\ref{eq21}).
\qed

\vskip2mm
 \noindent
{\bf Lemma 3.6.}  {\it $c_n$ is a critical value.}
 \vskip0mm
 \noindent
{\bf Proof.}\; By the boundedness of $\{(u_k,v_k)\}$ and the fact $X$ is a finite dimensional space, it follows that $\{(u_k,v_k)\}$ produces a convergent subsequence in $X$, which is renamed as $\{(u_k,v_k)\}$. Let $(u_n,v_n)=\lim\limits_{k\to\infty}(u_k,v_k)$. Note that $\varphi_V\in C^1(X,\R)$. Then Lemma 3.4 implies that
$$
\varphi_V'(u_n,v_n)=\lim\limits_{k\to\infty}\varphi_V'(u_k,v_k)=0\;\;\;\text{and}\;\;\;\varphi_V(u_n,v_n)=\lim\limits_{k\to\infty}\varphi_V(u_k,v_k)=c_n.
$$
Thus we finish the proof.
\qed

\vskip2mm
 \noindent
{\bf Proof of Theorem 3.1.}\; (i) By Lemma 3.6, we know that for each $n$ large enough, there is at least one solution $(u_n,v_n)$ of (\ref{eq1}) such that $\varphi_V(u_n,v_n)=c_n$, where $c_n$ is given by (\ref{eq20}). When $0<r_k\le R_n$, let $H_{r_k}=\{(c,d)\in\R\times\R:|(c,d)|=r_k\}+Q_V$ and $B_{r_k}=\{(u,v)\in \R^2:|(u,v)|\le r_k\}$. Let $\pi:X\to\R^2$ be the (continuous) projection of $X$ onto $\R^2$, which is defined by
$$
\pi(u,v)=\left(\frac{1}{|V|}\int_Vud\mu,\frac{1}{|V|}\int_Vvd\mu\right)\;\;\;\mbox{for}\;(u,v)\in X.
$$
Let $\gamma\in C(B_{R_n},X)$ such that $\gamma|_{\partial B_{R_n}}=id|_{\partial B_{R_n}}$. Define $\gamma_t:\R^2\to \R^2$ by
$$
\gamma_t(u,v)=t\pi(\gamma(u,v))+(1-t)(u,v),
$$
where $t\in [0,1]$. Thus $\gamma_0=id$, $\gamma_1=\pi\circ\gamma$ and $\gamma_t|_{\partial B_{R_n}}=id$ for all $t$. Next, we show that $\gamma(B_{R_n})\cap H_{r_k}\not=\emptyset$. Let $(u,v)$ be the element of the set $\{(c,d)\in\R\times\R:|(c,d)|=r_k\}$. Note that $r_k\le R_n$. Then $(u,v)\in B_{R_n}$. Hence, we have $deg(id,B_{R_n},(u,v))=1$. By homotopy invariance we have
$$
deg(\pi\circ\gamma,B_{R_n},(u,v))=deg(id,B_{R_n},(u,v))=1.
$$
Then there exists $(u^*,v^*)\in B_{R_n}$ such that
\begin{eqnarray}
\label{eq22}
\pi\circ\gamma(u^*,v^*)=(u,v)\in\partial B_{r_k}.
\end{eqnarray}
Suppose that $\gamma(u^*,v^*)=(u',v')\in X$. Then $(u',v')=(\tilde{u}'+\bar{u}',\tilde{v}'+\bar{v}')$, where $(\bar{u}',\bar{v}')\in \R^2$ and $(\tilde{u}',\tilde{v}')\in Q_V$, and we have $\pi(u',v')=(\bar{u}',\bar{v}')$. By (\ref{eq22}) we have $(u,v)=(\bar{u}',\bar{v}')$. Hence, $(\bar{u}',\bar{v}')\in \partial B_{r_k}$ and then $(u',v')=(\tilde{u}'+\bar{u}',\tilde{v}'+\bar{v}')\in Q_V+\partial  B_{r_k}$, that is $\gamma(u^*,v^*)\in Q_V+\partial B_{r_k}$. Then $\gamma$ intersects the hyperplane $H_{r_k}$ for any $\gamma\in S_n$. Then
$$
\max\limits_{(x,y)\in B_{R_n}}\varphi_V(\gamma(x,y))\geq\inf\limits_{(c+\tilde{u},d+\tilde{v})\in H_{r_k}}\varphi_V(c+\tilde{u},d+\tilde{v}).
$$
Therefore by using Lemma 3.3 we acquire that
$$
\lim\limits_{n\to\infty}c_n=+\infty,
$$
and (i) follows.
\par
(ii) For any fixed $m\in \mathbb N$, the subset $P_m$ of $X$ is defined by
$$
P_m=\{(u,v)\in X:(u,v)=(\bar{u},\bar{v})+(\tilde{u},\tilde{v}),\;(\bar{u},\bar{v})\in\R^2,\;|(\bar{u},\bar{v})|\leq r_m,\;(\tilde{u},\tilde{v})\in Q_V\}.
$$
By Lemma 3.1, we know that $\varphi_V$ is also coercive in $P_m\subset X$. Therefore $\varphi_V$ is bounded below on $P_m$.
We set
\begin{eqnarray}
\label{eq23}
\mu_m=\inf\limits_{(u,v)\in P_m}\varphi_V(u,v)
\end{eqnarray}
and assume that $\{(u_k,v_k)\}$ is a sequence in $P_m$ such that
\begin{eqnarray}
\label{eq24}
\varphi_V(u_k,v_k)\to\mu_m\;\;\;\mbox{as}\;k\to\infty.
\end{eqnarray}
Moreover, there are
$$
(u_k,v_k)=(\bar{u}_k,\bar{v}_k)+(\tilde{u}_k,\tilde{v}_k),\;\;\;(\bar{u}_k,\bar{v}_k)\in\R^2\;\;\;\text{and}\;|(\bar{u}_k,\bar{v}_k)|\le r_m.
$$
From the coerciveness of $\varphi_V$ in $P_m$ and (\ref{eq24}) we obtain that $\{(u_k,v_k)\}$ is a bounded sequence in $X$ and thus passing to a subsequence, which is renamed $\{(u_k,v_k)\}$, together with $dim X<\infty$, we have
$$
(u_k,v_k)\to (u_m^*,v_m^*)\;\;\text{in}\;\;X.
$$
Note that $P_m$ is a  closed subset of $X$. Then $(u_m^*,v_m^*)\in P_m$.
Furthermore, by the fact that $\varphi_V$ is continuous, we have
$$
\mu_m=\lim_{k\to\infty}\varphi_V(u_k,v_k)=\varphi_V(u_m^*,v_m^*).
$$
\par
If we can show $(u_m^*,v_m^*)$ is in the interior of $P_m$, then $(u_m^*,v_m^*)$ is a critical point of $\varphi_V$. To this end, let
$
(u_m^*,v_m^*)=(\bar{u}_m^*,\bar{v}_m^*)+(\tilde{u}_m^*,\tilde{v}_m^*)
$
with $(\bar{u}_m^*,\bar{v}_m^*)\in \R^2$ and $(\tilde{u}_m^*,\tilde{v}_m^*)\in Q_V$.
 If $R_n<r_m$, then $\{(a,b)\in\R^2:|(a,b)|=R_n\}\subset P_m$.  Then by the fact that $(u_m^*,v_m^*)\in P_m$, we have
$$
\mu_m=\varphi_V(u_m^*,v_m^*)\le\sup\limits_{(a,b)\in\R^2,|(a,b)|=R_n}\varphi_V(a,b)
$$
and then it follows from Lemma 3.2 that
$$
\lim\limits_{n\to\infty}\varphi_V(u_m^*,v_m^*)=-\infty,
$$
which, together with  Lemma 3.3, implies that  $|(\bar{u}_m^*,\bar{v}_m^*)|\not=r_m$ for large  $m$. So $(u_m^*,v_m^*)$ is in the interior of $P_m$. Hence, $(u_m^*,v_m^*)$ is a critical point of $\varphi_V$ and $\varphi_V(u_m^*,v_m^*)\to -\infty$ as $m\to\infty$. The proof is completed.
\qed

\vskip2mm
{\section{The locally finite graph case}}
  \setcounter{equation}{0}
\par
 In this section, we consider the existence of infinitely many solutions for system (\ref{eqA}) which is defined on locally finite graph. Being different from  (\ref{eq1}) which is defined on finite graph, we don't have to account for $h_i$, $i=1,2$ and the working space is changed, so that
 the assumptions on $F$ need to be modified correspondingly. We make the following assumptions:\\
  $(A)'$ \;  $F(x,t,s)$ is continuously differentiable in $(t,s)\in \R \times \R$ for all $x\in \Omega$;\\
$(H_3)$ \; $\int_\Omega F(x,0,0)d\mu=0$ and there exist functions $f_i:\Omega\to\R^+,i=1,2,3,4$ with $\frac{pq-q}{p}\|f_2\|_{\infty}+\|f_4\|_{\infty}<\frac{1}{C_{m_2,q(\Omega)}^q}$ and $\|f_1\|_{\infty}+\|f_2\|_{\infty}<\frac{1}{C_{m_1,p(\Omega)}^p},$ such that
\begin{eqnarray}
\begin{cases}
|F_u(x,u,v)|\le f_1(x)|u|^{p-1}+f_2(x)|v|^{\frac{pq-q}{p}},\\
|F_v(x,u,v)|\le f_3(x)|u|^{\frac{pq-p}{q}}+f_4(x)|v|^{q-1},\\
\end{cases}
\end{eqnarray}
for all $(x,u,v)\in \Omega\times\R\times\R$, where $C_{m_1,p(\Omega)}=\frac{C_p}{\mu_{\min,\Omega}}(1+|\sum_{x\in\Omega}\mu(x)|)$ with $C_p$ satisfying $\|u\|_{L^{\theta}(\Omega)}\le C_p\|u\|_{W_0^{m_1,p}(\Omega)}$, $C_{m_2,q(\Omega)}=\frac{C_q}{\mu_{\min,\Omega}}(1+|\sum_{x\in\Omega}\mu(x)|)$ with $C_q$ satisfying $\|u\|_{L^{\theta}(\Omega)}\le C_q\|u\|_{W_0^{m_2,q}(\Omega)}$;\\
$(H_4)$ \; there hold
\begin{eqnarray}
\label{eqA3}
\limsup\limits_{R\rightarrow +\infty}\inf\limits_{a,b\in \R,|(a,b)|=R}\int_{\Omega}F(x,a,b)d\mu=+\infty
\end{eqnarray}
and
\begin{eqnarray}
\label{eqA4}
& &\liminf\limits_{r\rightarrow +\infty}\sup\limits_{c,d\in \R,|(c,d)|=r}\frac{\int_\Omega F(x,c,d)d\mu}{|c|^p+|d|^q}<\min\Bigg\{\int_\Omega\left(\frac{1}{p}f_1(x)-\frac{1}{p}f_1(x)2^{p-1}-\frac{q-1}{q}f_3(x)\right)d\mu,\notag\\
& &\int_\Omega\left(\frac{1}{q}f_4(x)-\frac{p-1}{p}f_2(x)2^{q-1}-\frac{1}{q}f_4(x)2^{q-1}\right)d\mu\Bigg\}.
\end{eqnarray}

\vskip2mm
\par
 We shall work in the space $W:=W_0^{m_1,p}(\Omega)\times W_0^{m_2,q}(\Omega)$ with the norm $\|(u,v)\|=\|u\|_{W_0^{m_1,p}(\Omega)}+\|v\|_{W_0^{m_2,q}(\Omega)}$. Then $W$ is a Banach space with finite dimension. The functional $\varphi_\Omega:W\to\R$ is defined as
\begin{eqnarray}
\label{eqA15}
\varphi_\Omega(u,v)=\frac{1}{p}\int_{\Omega\cup\partial\Omega}|\nabla^{m_1}u|^pd\mu+\frac{1}{q}\int_{\Omega\cup\partial\Omega}|\nabla^{m_2}v|^qd\mu-\int_{\Omega} F(x,u,v)d\mu.
\end{eqnarray}
Then the assumption $(A)$ implies that $\varphi_{\Omega}\in C^1(W,\R)$ and
\begin{eqnarray}
\label{eqA16}
\langle\varphi_\Omega'(u,v),(\phi_1,\phi_2)\rangle&=&\int_{\Omega\cup\partial\Omega}\left[(\pounds_{m_1,p}u)\phi_1-F_u(x,u,v)\phi_1\right]d\mu\nonumber\\
&&+\int_{\Omega\cup\partial\Omega}\left[(\pounds_{m_2,q}v)\phi_2-F_v(x,u,v)\phi_2\right]d\mu
\end{eqnarray}
for any $(u,v),(\phi_1,\phi_2)\in W$. The problem about finding the solutions for system (\ref{eqA}) is attributed to finding the critical points of the functional $\varphi_{\Omega}$ on $W$ (see \cite{Yu 2023}). The following is our main result in this section.
\vskip2mm
\noindent
{\bf Theorem 4.1.} {\it Suppose $(A)'$, $(H_3)$ and  $(H_4)$ hold. Then\\
$(i)$ \; system (\ref{eqA}) possesses a sequence of solutions  $\{(u_n,v_n)\}$ such that $(u_n,v_n)$ is a critical point of $\varphi_{\Omega}$ and $\varphi_{\Omega}(u_n,v_n)=+\infty$  as $n\to \infty$;\\
$(ii)$ \; system (\ref{eqA}) possesses a sequence of $\{(u_m^*,v_m^*)\}$ such that $(u_m^*,v_m^*)$ is a local minimum point of $\varphi_{\Omega}$ and $\varphi_{\Omega}(u_m^*,m_m^*)=-\infty$ as $m\to \infty$.}

\vskip2mm
\par
The line of proofs for Theorem 4.1 is similar to Theorem 3.1 except for some details caused by the different norms in Lemma 4.1-Lemma 4.3 below. Consider the direct sum decomposition $W=\R^2\oplus Q_{\Omega}$ with $Q_{\Omega}=\{(u,v)\in W:\int_{\Omega} u(x)d\mu=0,\int_{\Omega}v(x) d\mu=0\}$. So for all $(u,v)\in W$ we can write $(u,v)=(\bar{u},\bar{v})+(\tilde{u},\tilde{v})$ with $(\bar{u},\bar{v})\in \R^2$ and $(\tilde{u},\tilde{v})\in Q_{\Omega}$.

 \vskip2mm
 \noindent
{\bf Lemma  4.1.}  {\it $\varphi_{\Omega}$ is coercive on $W$.}
 \vskip0mm
 \noindent
{\bf Proof.}\;Note that $\int_{\Omega}F(x,0,0)d\mu=0$. Similar to the proof in Lemma 3.1, we have
\begin{eqnarray*}
\int_{\Omega}F(x,s,t)d\mu&\le&\frac{1}{p}\|f_1\|_{\infty}\|s\|_p^p+\|f_2\|_{\infty}\int_{\Omega}|t|^{\frac{pq-q}{p}}|s|d\mu+\frac{1}{q}\|f_4\|_{\infty}\|t\|_q^q.
\end{eqnarray*}
Hence, together with Young's inequality, we have
\begin{eqnarray*}
\varphi_{\Omega}(u,v)&\ge&\frac{1}{p}\|u\|_{W_0^{m_1,p}(\Omega)}^p+\frac{1}{q}\|v\|_{W_0^{m_2,q}(\Omega)}^q-\frac{1}{p}\|f_1\|_{\infty}\|u\|_p^p-\frac{p-1}{p}\|f_2\|_{\infty}\|v\|_q^q-\frac{1}{p}\|f_2\|_{\infty}\|u\|_p^p-\frac{1}{q}\|f_4\|_{\infty}\|v\|_q^q\\
&\ge&\frac{1}{p}\|u\|_{W_0^{m_1,p}(\Omega)}^p+\frac{1}{q}\|v\|_{W_0^{m_2,q}(\Omega)}^q-\frac{1}{p}\|f_1\|_{\infty}C_{m_1,p(\Omega)}^p\|u\|_{W_0^{m_1,p}(\Omega)}^p-\frac{p-1}{p}\|f_2\|_{\infty}C_{m_2,q(\Omega)}^q\|v\|_{W_0^{m_2,q}(\Omega)}^q\\
& &-\frac{1}{p}\|f_2\|_{\infty}C_{m_1,p(\Omega)}^p\|u\|_{W_0^{m_1,p}(\Omega)}^p-\frac{1}{q}\|f_4\|_{\infty}C_{m_2,q(\Omega)}^q\|v\|_{W_0^{m_2,q}(\Omega)}^q.
\end{eqnarray*}
Since $\frac{pq-q}{p}\|f_2\|_{\infty}+\|f_4\|_{\infty}<\frac{1}{C_{m_2,q(\Omega)}^q}$, $\|f_1\|_{\infty}+\|f_2\|_{\infty}<\frac{1}{C_{m_1,p(\Omega)}^p}$ and $p,q>1$, then $\varphi_{\Omega}(u,v)$ is coercive on $W$.
\qed

 \vskip2mm
 \noindent
{\bf Lemma 4.2.}  {\it There is a sequence $\{R_n\}$ such that
$$
\lim\limits_{n\to\infty}R_n=+\infty,\;\lim\limits_{n\to\infty}\left[\sup\limits_{a,b\in\R,|(a,b)|=R_n}\varphi_{\Omega}(a,b)\right]=-\infty.
$$}
 \vskip0mm
 \noindent
{\bf Proof.}\; Note that $|\nabla^{m_1} a|=|\nabla^{m_2} b|=0$ for all $a,b\in\R$. Since
$$
\varphi_{\Omega}(a,b)=-\int_{\Omega} F(x,a,b) d\mu,
$$
then  the proof is completed by  (\ref{eqA3}).\qed

\vskip2mm
 \noindent
{\bf Lemma 4.3.}  {\it There is a sequence $\{r_m\}$ such that
$$
\lim\limits_{m\to\infty}r_m=+\infty,\;\lim\limits_{m\to\infty}\left[\inf\limits_{c,d\in\R,|(c,d)|=r_m,(\tilde{u},\tilde{v})\in Q_{\Omega}}\varphi_{\Omega}(c+\tilde{u},d+\tilde{v})\right]=+\infty.
$$}
 \vskip0mm
 \noindent
{\bf Proof.}\; For any given $c,d\in\R$, $|(c,d)|=r_m$, and $(\tilde{u},\tilde{v})\in Q_{\Omega}$, set $(u,v)=(c+\tilde{u},d+\tilde{v})$. Similar to the proof in Lemma 3.3,  we have
\begin{eqnarray*}
& &\varphi_{\Omega}(u,v)\\
&=&\varphi_{\Omega}(c+\tilde{u},d+\tilde{v})\\
&=&\frac{1}{p}\int_{\Omega\cup\partial\Omega}|\nabla^{m_1}(c+\tilde{u})|^pd\mu+\frac{1}{q}\int_{\Omega\cup\partial\Omega}|\nabla^{m_2}(d+\tilde{v})|^q d\mu-\int_{\Omega} F(x,c+\tilde{u},d+\tilde{v})d\mu\\
&\ge&\frac{1}{p}\int_{\Omega\cup\partial\Omega}|\nabla^{m_1}\tilde{u}|^pd\mu+\frac{1}{q}\int_{\Omega\cup\partial\Omega}|\nabla^{m_2}\tilde{v}|^q d\mu-\int_{\Omega}F(x,c,d)d\mu\\
& &-\int_{\Omega}\bigg[\frac{1}{p}(f_1(x)2^{p-1}+f_2(x))|\tilde{u}|^p+\frac{1}{q}\left(f_3(x)+2^{q-1}f_4(x)+2^{q-1}\frac{pq-q}{p}f_2(x)\right)|\tilde{v}|^q\bigg]d\mu\\
& &-|c|^p\int_{\Omega}\left(\frac{1}{p}f_1(x)2^{p-1}-\frac{1}{p}f_1(x)+\frac{q-1}{q}f_3(x)\right)d\mu-|d|^q\int_{\Omega}\left(\frac{p-1}{p}f_2(x)2^{q-1}+\frac{1}{q}f_4(x)2^{q-1}-\frac{1}{q}f_4(x)\right)d\mu\\
&\ge&\frac{1}{p}\int_{\Omega\cup\partial\Omega}|\nabla^{m_1}\tilde{u}|^p d\mu+\frac{1}{q}\int_{\Omega\cup\partial\Omega}|\nabla^{m_2}\tilde{v}|^q d\mu-\int_{\Omega}F(x,c,d)d\mu-\frac{1}{p}(\|f_1\|_{\infty}2^{p-1}+\|f_2\|_{\infty})C_{m_1,p(\Omega)}^p\|\tilde{u}\|_{W_0^{m_1,p}(\Omega)}^p\\
& &-\frac{1}{q}\left(
\|f_3\|_{\infty}+2^{q-1}\|f_4\|_{\infty}+2^{q-1}\frac{pq-q}{p}\|f_2\|_{\infty}\right)C_{m_2,q(\Omega)}^q\|\tilde{v}\|_{W_0^{m_2,q}(\Omega)}^q\\
& &-|c|^p\int_{\Omega}\left(\frac{1}{p}f_1(x)2^{p-1}-\frac{1}{p}f_1(x)+\frac{q-1}{q}f_3(x)\right)d\mu-|d|^q\int_{\Omega}\left(\frac{p-1}{p}f_2(x)2^{q-1}+\frac{1}{q}f_4(x)2^{q-1}-\frac{1}{q}f_4(x)\right)d\mu\\
&=&\frac{1}{p}\left[1-(\|f_1\|_{\infty}2^{p-1}+\|f_2\|_{\infty})C_{m_1,p(\Omega)}^p\right]\int_{\Omega\cup\partial\Omega}|\nabla^{m_1}\tilde{u}|^p d\mu\\
& &+\frac{1}{q}\bigg[1-\left(\|f_3\|_{\infty}+2^{q-1}\|f_4\|_{\infty}+2^{q-1}\frac{pq-q}{p}\|f_2\|_{\infty}\right)C_{m_2,q(\Omega)}^q\bigg]\int_{\Omega\cup\partial\Omega}|\nabla^{m_2}\tilde{v}|^q d\mu\\
& &-|c|^p\int_{\Omega}\left(\frac{1}{p}f_1(x)2^{p-1}-\frac{1}{p}f_1(x)+\frac{q-1}{q}f_3(x)\right)d\mu\\
& &-|d|^q\int_{\Omega}\left(\frac{p-1}{p}f_2(x)2^{q-1}+\frac{1}{q}f_4(x)2^{q-1}-\frac{1}{q}f_4(x)\right)d\mu-\int_{\Omega} F(x,c,d)d\mu\\
&\ge&\frac{1}{p}\left[1-(\|f_1\|_{\infty}2^{p-1}+\|f_2\|_{\infty})C_{m_1,p(\Omega)}^p\right]\int_{\Omega\cup\partial\Omega}|\nabla^{m_1}\tilde{u}|^p d\mu\\
& &+\frac{1}{q}\bigg[1-\left(\|f_3\|_{\infty}+2^{q-1}\|f_4\|_{\infty}+2^{q-1}\frac{pq-q}{p}\|f_2\|_{\infty}\right)C_{m_2,q(\Omega)}^q\bigg]\int_{\Omega\cup\partial\Omega}|\nabla^{m_2}\tilde{v}|^q d\mu\\
& &+(|c|^p+|d|^q)\Bigg(\min\bigg\{\int_{\Omega}\left(\frac{1}{p}f_1(x)-\frac{1}{p}f_1(x)2^{p-1}-\frac{q-1}{q}f_3(x)\right)d\mu,\\
& &\int_{\Omega}\left(\frac{1}{q}f_4(x)-\frac{p-1}{p}f_2(x)2^{q-1}-\frac{1}{q}f_4(x)2^{q-1}\right)d\mu\bigg\}-\frac{\int_{\Omega}F(x,c,d)d\mu}{|c|^p+|d|^q}\Bigg).
\end{eqnarray*}
Thus, the proof is completed by  (\ref{eqA4}).\qed
\par
 \vskip2mm
Consider the set $S_n=\{\gamma\in C(B_{R_n},W),\gamma|_{\partial B_{R_n}}=id|_{\partial B_{R_n}}\}$, where $R_n$ is defined in Lemma 4.2 and $B_{R_n}=\{(u,v)\in \R^2:|(u,v)|\le R_n\}$. Define
\begin{eqnarray}
\label{eqA20}
c_n=\inf_{\gamma\in S_n}\left[\max_{(x,y)\in B_{R_n}}\varphi_{\Omega}(\gamma(x,y))\right].
\end{eqnarray}
Let $\pi_\Omega:X\to \R^2$ be the (continuous) projection of $W$ onto $\R^2$, which is defined by
$$
\pi_\Omega(u,v)=\left(\frac{1}{|\Omega|}\int_{\Omega} u d\mu,\frac{1}{|\Omega|}\int_{\Omega} v d\mu\right)\;\;\;\mbox{for}\;(u,v)\in W.
$$
The rest of proofs of Theorem 4.1 is the same as Theorem 3.1 with replacing $X$, $Q_V$, $\phi_V$ and $\pi_V$ with $W$, $Q_\Omega$, $\phi_\Omega$ and  $\pi_\Omega$, respectively.
\qed

\vskip2mm
{\section{Examples}}
  \setcounter{equation}{0}
 \vskip0mm
 \noindent
{\bf Example 5.1.}\; Let $p=2$ and $q=3$. Consider the following system:
\begin{eqnarray}
\begin{cases}
\label{eq51}
  \pounds_{m_1,2}u+h_1(x)u=F_u(x,u,v),\;\;\;\;\hfill x\in V,\\
  \pounds_{m_2,3}v+h_2(x)v=F_v(x,u,v),\;\;\;\;\hfill x\in V,\\
\end{cases}
\end{eqnarray}
where $G=(V,E)$ is a finite graph, $m_i\ge 1,i=1,2$ are integers, the measure $\mu(x)\ge\mu_0=\frac{1}{4000}$ for all $x\in V$, $h_1(x)=h_2(x)=546$ and $\int_Vd\mu=\frac{1}{300}$,
$$
F(x,u,v)=(|u|^2+|v|^3)\sin(\ln(|u|^2+|v|^3+1)).
$$
Next, we verify that $F$ satisfies the conditions in Theorem 3.1:
\par
\begin{itemize}
\item Let
$$
f_1(x)\equiv 4,f_2(x)\equiv f_3(x)\equiv 1,f_4(x)\equiv 6.
$$
Then
\begin{eqnarray*}
|F_u(x,u,v)|&=&\left|2|u|\sin(\ln(|u|^2+|v|^3+1))+(|u|^2+|v|^3)\cos(\ln(|u|^2+|v|^3+1))\frac{2|u|}{|u|^2+|v|^3+1}\right|\\
&\le&2|u||\sin(\ln(|u|^2+|v|^3+1))|+\frac{2|u|(|u|^2+|v|^3)}{|u|^2+|v|^3+1}|\cos(\ln(|u|^2+|v|^3+1))|\\
&\le&2|u|+\frac{2|u|(|u|^2+|v|^3+1)}{|u|^2+|v|^3+1}\\
&=&4|u|\\
&\le&4|u|+|v|^2
\end{eqnarray*}
and
\begin{eqnarray*}
|F_v(x,u,v)|&=&\left|3|v|^2\sin(\ln(|u|^2+|v|^3+1))+(|u|^2+|v|^3)\cos(\ln(|u|^2+|v|^3+1))\frac{3|v|^2}{|u|^2+|v|^3+1}\right|\\
&\le&3|v|^2|\sin(\ln(|u|^2+|v|^3+1))|+\frac{3|v|^2(|u|^2+|v|^3)}{|u|^2+|v|^3+1}|\cos(\ln(|u|^2+|v|^3+1))|\\
&\le&3|v|^2+\frac{3|v|^2(|u|^2+|v|^3+1)}{|u|^2+|v|^3+1}\\
&=&6|v|^2\\
&\le&|u|+6|v|^2.
\end{eqnarray*}
\item $\int_V F(x,0,0)d\mu=0$.
\item We have
\begin{eqnarray*}
\frac{pq-q}{p}\|f_2\|_{\infty}+\|f_4\|_{\infty}=7.5<\frac{1}{K_q^q}=\frac{1}{K_3^3}=\frac{(\mu_{min}h_{min})^{\frac{3}{2}}}{\sum_{x\in V}\mu(x)}\approx15.1,\\
\|f_1\|_{\infty}+\|f_2\|_{\infty}=5<\frac{1}{K_p^p}=\frac{1}{K_2^2}=\frac{\mu_{min}h_{min}}{(\sum_{x\in V}\mu(x))^{\frac{2}{3}}}\approx6.1.
\end{eqnarray*}
\item
$$\frac{\int_V F(x,a,b)d\mu}{|a|^2+|b|^3}=\frac{\int_V (|a|^2+|b|^3)\sin(\ln(|a|^2+|b|^3+1))d\mu}{|a|^2+|b|^3}=\sin(\ln(|a|^2+|b|^3+1)).$$
Let
$$R_k=\sqrt\frac{e^{2k\pi+1/2\pi}-1}{2},R_k'=\sqrt[3]{\frac{e^{2k\pi+1/2\pi}-1}{2}},r_k=\sqrt\frac{e^{2k\pi+3/2\pi}-1}{2},r_k'=\sqrt[3]{\frac{e^{2k\pi+3/2\pi}-1}{2}}$$
for $k=1,2,...$ Then
$$
\lim\limits_{k\to+\infty}\frac{\int_V F(x,R_k,R_k')d\mu}{|R_k|^2+|R_k'|^3}=\lim\limits_{k\to+\infty}\sin(\ln(|R_k|^2+|R_k'|^3+1))=\lim\limits_{k\to+\infty}1=1
$$
and
$$
\lim\limits_{k\to+\infty}\frac{\int_V F(x,r_k,r_k')d\mu}{|r_k|^2+|r_k'|^3}=\lim\limits_{k\to+\infty}\sin(\ln(|r_k|^2+|r_k'|^3+1))=\lim\limits_{k\to+\infty}-1=-1.
$$
We also have
\begin{eqnarray*}
\max\left\{\frac{1}{p}\int_V h_1(x)d\mu,\frac{1}{q}\int_V h_2(x)d\mu\right\}=0.91
\end{eqnarray*}
and
\begin{eqnarray*}
& &\min\Bigg\{\int_V\left(\frac{1}{p}f_1(x)-\frac{1}{p}f_1(x)2^{p-1}-\frac{q-1}{q}f_3(x)-\frac{1}{p}h_1(x)\right)d\mu,\\
& &\int_V\left(\frac{1}{q}f_4(x)-\frac{p-1}{p}f_2(x)2^{q-1}-\frac{1}{q}f_4(x)2^{q-1}-\frac{1}{q}h_2(x)\right)d\mu\Bigg\}\approx-0.92.
\end{eqnarray*}
As a result, there hold
\begin{eqnarray*}
\limsup\limits_{R\rightarrow +\infty}\inf\limits_{a,b\in \R,|(a,b)|=R}\frac{\int_V F(x,a,b)d\mu}{|a|^p+|b|^q} >\max\left\{\frac{1}{p}\int_Vh_1(x)d\mu,\frac{1}{q}\int_Vh_2(x)d\mu\right\}
\end{eqnarray*}
and
\begin{eqnarray*}
& &\liminf\limits_{r\rightarrow +\infty}\sup\limits_{c,d\in \R,|(c,d)|=r}\frac{\int_V F(x,c,d)d\mu}{|c|^p+|d|^q}<\min\Bigg\{\int_V\left(\frac{1}{p}f_1(x)-\frac{1}{p}f_1(x)2^{p-1}-\frac{q-1}{q}f_3(x)-\frac{1}{p}h_1(x)\right)d\mu,\notag\\
& &\int_V\left(\frac{1}{q}f_4(x)-\frac{p-1}{p}f_2(x)2^{q-1}-\frac{1}{q}f_4(x)2^{q-1}-\frac{1}{q}h_2(x)\right)d\mu\Bigg\}.
\end{eqnarray*}

\end{itemize}
Hence, the conditions in Theorem 3.1 are satisfied and we have the following results:\\
$(i)$ \; system  (\ref{eq51}) possesses a sequence of solutions  $\{(u_n,v_n)\}$ such that $(u_n,v_n)$ is a mountain pass type critical point of $\varphi_V$ and $\varphi_V(u_n,v_n)=+\infty$  as $n\to \infty$;\\
$(ii)$ \; system  (\ref{eq51}) possesses  a sequence of $\{(u_m^*,v_m^*)\}$  such that $(u_m^*,v_m^*)$ is a local minimum point of $\varphi_V$ and $\varphi_V(u_m^*,v_m^*)=-\infty$ as $m\to \infty$.
\par

\vskip2mm
 \noindent
{\bf Example 5.2.}\; Let $m_1=m_2=1$, $p=2$ and $q=3$. Consider the following system:
\begin{eqnarray}
\begin{cases}
\label{eq53}
  \pounds_{1,2}u=F_u(x,u,v),\;\;\;\;\hfill x\in \Omega,\\
  \pounds_{1,3}v=F_v(x,u,v),\;\;\;\;\hfill x\in \Omega,\\
\end{cases}
\end{eqnarray}
where $G=(V,E)$ is a locally finite graph, $\Omega$ is a bounded domain of $V$ such that $\Omega^{\circ}\neq\emptyset$ and there is at least one $y\in \partial \Omega$ satisfying that $y \thicksim x$ for each $x\in\Omega$, the measure $\mu(x)\ge\mu_0=\frac{1}{30000}$ for all $x\in \Omega$, $\mu(x)\le\mu_{\max}=\frac{1}{6000}$ for all $x\in \Omega\cup\partial \Omega$, the weight $w_{xy}\ge w_{\min}=51$ for all $x\in \Omega\cup\partial \Omega$ and $\int_{\Omega}d\mu=\frac{1}{300}$,
$$
F(x,u,v)=(|u|^2+|v|^3)\sin(\ln(|u|^2+|v|^3+1)).
$$
Next, we verify that $F$ satisfies the conditions in Theorem 4.1:
\par
\begin{itemize}
\item Let
$$
f_1(x)\equiv 4,f_2(x)\equiv f_3(x)\equiv 1,f_4(x)\equiv 6.
$$
Then
\begin{eqnarray*}
|F_u(x,u,v)|&=&\left|2|u|\sin(\ln(|u|^2+|v|^3+1))+(|u|^2+|v|^3)\cos(\ln(|u|^2+|v|^3+1))\frac{2|u|}{|u|^2+|v|^3+1}\right|\\
&\le&2|u||\sin(\ln(|u|^2+|v|^3+1))|+\frac{2|u|(|u|^2+|v|^3)}{|u|^2+|v|^3+1}|\cos(\ln(|u|^2+|v|^3+1))|\\
&\le&2|u|+\frac{2|u|(|u|^2+|v|^3+1)}{|u|^2+|v|^3+1}\\
&=&4|u|\\
&\le&4|u|+|v|^2
\end{eqnarray*}
and
\begin{eqnarray*}
|F_v(x,u,v)|&=&\left|3|v|^2\sin(\ln(|u|^2+|v|^3+1))+(|u|^2+|v|^3)\cos(\ln(|u|^2+|v|^3+1))\frac{3|v|^2}{|u|^2+|v|^3+1}\right|\\
&\le&3|v|^2|\sin(\ln(|u|^2+|v|^3+1))|+\frac{3|v|^2(|u|^2+|v|^3)}{|u|^2+|v|^3+1}|\cos(\ln(|u|^2+|v|^3+1))|\\
&\le&3|v|^2+\frac{3|v|^2(|u|^2+|v|^3+1)}{|u|^2+|v|^3+1}\\
&=&6|v|^2\\
&\le&|u|+6|v|^2.
\end{eqnarray*}
\item $\int_{\Omega} F(x,0,0)d\mu=0$.
\item Applying Lemma 2.2 in \cite{Yang 2023}, we have
\begin{eqnarray*}
\frac{pq-q}{p}\|f_2\|_{\infty}+\|f_4\|_{\infty}=5<\frac{1}{C_{m_2,q(\Omega)}^q}=\frac{\mu_0w_{\min}^{\frac{3}{2}}}{(2\mu_{\max})^{\frac{3}{2}}(1+|\sum_{x\in \Omega}\mu(x)|)^3}\approx1975.06,\\
\|f_1\|_{\infty}+\|f_2\|_{\infty}=5<\frac{1}{C_{m_1,p(\Omega)}^p}=\frac{\mu_0w_{\min}}{2\mu_{\max}(1+|\sum_{x\in \Omega}\mu(x)|)^2}\approx5.07.
\end{eqnarray*}
\item
$$\frac{\int_{\Omega} F(x,c,d)d\mu}{|c|^2+|d|^3}=\frac{\int_{\Omega}(|c|^2+|d|^3)\sin(\ln(|c|^2+|d|^3+1))d\mu}{|c|^2+|d|^3}=\sin(\ln(|c|^2+|d|^3+1)).$$
Let
$$R_k=\sqrt\frac{e^{2k\pi+1/2\pi}-1}{2},R_k'=\sqrt[3]{\frac{e^{2k\pi+1/2\pi}-1}{2}},r_k=\sqrt\frac{e^{2k\pi+3/2\pi}-1}{2},r_k'=\sqrt[3]{\frac{e^{2k\pi+3/2\pi}-1}{2}}$$
for $k=1,2,...$ Then
$$
\lim\limits_{k\to+\infty}\int_{\Omega} F(x,R_k,R_k')d\mu=\lim\limits_{k\to+\infty}\frac{1}{300}(|R_k|^2+|R_k'|^3)\sin(\ln(|R_k|^2+|R_k'|^3+1))=+\infty
$$
and
$$
\lim\limits_{k\to+\infty}\frac{\int_{\Omega} F(x,r_k,r_k')d\mu}{|r_k|^2+|r_k'|^3}=\lim\limits_{k\to+\infty}\sin(\ln(|r_k|^2+|r_k'|^3+1))=\lim\limits_{k\to+\infty}-1=-1.
$$
We also have
\begin{eqnarray*}
& &\min\Bigg\{\int_{\Omega}\left(\frac{1}{p}f_1(x)-\frac{1}{p}f_1(x)2^{p-1}-\frac{q-1}{q}f_3(x)\right)d\mu,\\
& &\int_{\Omega}\left(\frac{1}{q}f_4(x)-\frac{p-1}{p}f_2(x)2^{q-1}-\frac{1}{q}f_4(x)2^{q-1}\right)d\mu\Bigg\}\approx-0.03.
\end{eqnarray*}
As a result, there hold
\begin{eqnarray*}
\limsup\limits_{R\rightarrow +\infty}\inf\limits_{a,b\in \R,|(a,b)|=R}\int_{\Omega} F(x,a,b)d\mu=+\infty
\end{eqnarray*}
and
\begin{eqnarray*}
& &\liminf\limits_{r\rightarrow +\infty}\sup\limits_{c,d\in \R,|(c,d)|=r}\frac{\int_{\Omega} F(x,c,d)d\mu}{|c|^p+|d|^q}<\min\Bigg\{\int_{\Omega}\left(\frac{1}{p}f_1(x)-\frac{1}{p}f_1(x)2^{p-1}-\frac{q-1}{q}f_3(x)\right)d\mu,\notag\\
& &\int_{\Omega}\left(\frac{1}{q}f_4(x)-\frac{p-1}{p}f_2(x)2^{q-1}-\frac{1}{q}f_4(x)2^{q-1}\right)d\mu\Bigg\}.
\end{eqnarray*}

\end{itemize}
Hence, the conditions in Theorem 4.1 are satisfied and we have the following results:\\
$(i)$ \; system  (\ref{eq53}) possesses a sequence of solutions  $\{(u_n,v_n)\}$ such that $(u_n,v_n)$ is a mountain pass type critical point of $\varphi_{\Omega}$ and $\varphi_{\Omega}(u_n,v_n)=+\infty$  as $n\to \infty$;\\
$(ii)$ \; system  (\ref{eq53}) possesses  a sequence of $\{(u_m^*,v_m^*)\}$  such that $(u_m^*,v_m^*)$ is a local minimum point of $\varphi_{\Omega}$ and $\varphi_{\Omega}(u_m^*,v_m^*)=-\infty$ as $m\to \infty$.

\vskip2mm
{\section{Remark}}
  \setcounter{equation}{0}
The conclusion $(ii)$ in Theorem 3.1 and Theorem 4.1 implies that $\varphi$ has no ground state energy and so the system (\ref{eq1}) and (\ref{eqA})  has no ground state solution under our assumptions. Moreover, in \cite{Zhang 2022}, Zhang-Zhang-Xie-Yu obtained that the system (\ref{eq1}) has at least $\mbox{dim} W$  solutions, where $W$ is the working space,  when $F$ has super-$(p,q)$-linearity and symmetry by using the symmetric mountain pass theorem. They could not obtain that the system has infinitely many solutions because of the nature of symmetric mountain pass theorem and the finiteness of dimension of $W$. Our results obtain the system (\ref{eq1})  has infinitely many solutions when $F$ has sub-$(p,q)$-linearity and no symmetry. To the best our knowledge, it seems that there are no results on infinitely many solutions for  the system (\ref{eq1}) on graph except for \cite {Pang 2024} where the results are different from our  Theorem 3.1 and Theorem 4.1.

 \vskip2mm
 \noindent
 {\bf Acknowledgments}\\
This work is supported by Yunnan Fundamental Research Projects of China (grant No: 202301AT070465) and  Xingdian Talent
Support Program for Young Talents of Yunnan Province in China.

\vskip2mm
\renewcommand\refname{References}
{}

\end{document}